\documentclass [12pt]{amsart}

\usepackage{amssymb,amsxtra,amsfonts}
\usepackage{epsf}           %
\usepackage{epsfig}             %
\usepackage{graphics}

\usepackage[colorlinks]{hyperref}

\usepackage{verbatim}  
\usepackage{bbding}   

\newenvironment{ppb}[1]
{\ \!\!\!\!\!\!\!\!\!\!\!\!\!\!\!\!\!\!\!\!\!\!\!\!\!\!\!\!\!\!\!\!\!\!\!\!\!\!\!\! {\bf PPB------------------------------------------------------------------------------------------------PPB}\newline \tiny {#1}
\  \newline\normalsize\phantom{f}\!\!\!\!\!\!\!\!\!\!\!\!\!\!\!\!\!\!\!\!\!\!\!\!\!\!\!\!\!\!\!\!\!\!\!\!\!\!\!\! {\bf PPB------------------------------------------------------------------------------------------------PPB}\newline}{}

\def\reE@DeclareMathSymbol#1#2#3#4{%
    \let#1=\undefined
    \DeclareMathSymbol{#1}{#2}{#3}{#4}}
\DeclareSymbolFont{symbolsC}{U}{txsyc}{m}{n}
\SetSymbolFont{symbolsC}{bold}{U}{txsyc}{bx}{n}
\DeclareFontSubstitution{U}{txsyc}{m}{n}
\reE@DeclareMathSymbol{\strictiff}{\mathrel}{symbolsC}{76}

\newcommand\beq{\begin{equation}}
\newcommand\eeq{\end{equation}}
\newcommand\bal{\begin{align*}}
\newcommand\eal{\end{align*}}   
\newcommand\bmx{\left(\begin{matrix}}
\newcommand\emx{\end{matrix}\right)}
\newcommand\bsmx{\left(\begin{smallmatrix}}
\newcommand\esmx{\end{smallmatrix}\right)}

\newcommand{\st}{\ \bigl\vert\ }

\def\part#1{\frac{\partial\phantom{q}}{\partial#1}}

\newcommand{\union}{\cup} 
 
\newcommand{\sdp}{{\ltimes}}

\newcommand{\MB}{\mathcal{M}_{\text{\rm B}}}

\newcommand{\Lie}{{\mathop{\rm Lie}}}

\newcommand{\Sym}{\mathop{\rm Sym}}

\DeclareMathOperator{\pr}{pr}

\newcommand{\Prod}{\prod}

\newcommand{\tr}{{\mathop{\rm Tr}}}
\newcommand{\Tr}{{\mathop{\rm Tr}}}
\DeclareMathOperator{\Hom}{Hom}         
\DeclareMathOperator{\Aut}{\mathop{\rm Aut}}

\newcommand{\SL}{{\mathop{\rm SL}}}
\newcommand{\PSL}{{\mathop{\rm PSL}}}
\newcommand{\GL}{{\mathop{\rm GL}}}

\newcommand{\SO}{{\mathop{\rm SO}}}

\newcommand{\Spin}{\mathop{\rm Spin}}
\newcommand{\PSU}{\mathop{\rm PSU}}

\renewcommand{\Im}{\mathop{\rm Im}}

\newcommand{\diag}{{\mathop{\rm diag}}}

\newcommand{\hk}{{hyperk\"ahler }}   

\newcommand{\IC}{\mathbb{C}}

\newcommand{\IF}{\mathbb{F}}

\newcommand{\IM}{\mathbb{M}}

\newcommand{\IO}{\mathbb{O}}
\newcommand{\IP}{\mathbb{P}}                                     
                           
\newcommand{\IR}{\mathbb{R}}

\newcommand{\IZ}{\mathbb{Z}}

\newcommand{\cC}{\mathcal{C}}

\newcommand{\cD}{\mathcal{D}}

\newcommand{\cM}{\mathcal{M}}

\newcommand{\cO}{\mathcal{O}}

\newcommand{\bcC}{\boldsymbol{\mathcal{C}}}

\newcommand{\g}{       \mathfrak{g}     }

\newcommand{\lt}{\mathfrak{t}}

\newcommand{\wh}{\widehat}

\newcommand{\al}{\alpha}

\newcommand{\be}{\beta}

\newcommand{\Si}{\Sigma}
\renewcommand{\th}{\theta}

\renewcommand{\bar}{\overline}

\makeatletter
 \newlength{\typesize}
\makeatother

\newlength{\vvoff}
\newlength{\hhoff}

\newcommand{\pf}{\begin{bpf}}

\newcommand{\pfms}{\begin{bpfms}}
\newcommand{\epf}{\end{bpf}\hfill$\square$\\}           
\newcommand{\epfms}{\end{bpfms}\hfill$\square$\\}       

\newcommand{\idea}{\begin{bidea}}

\newcommand{\eidea}{\end{bidea}\hfill$\square$\\}           

\newcommand{\sk}{\begin{bsk}}    

\newcommand{\esk}{\end{bsk}\hfill$\square$\\}           
\newcommand{\sketch}{\begin{bsketch}}

\newcommand{\esketch}{\end{bsketch}\hfill$\square$\\}

\newtheorem {hypo}{\bf\hspace{-\parindent}Hypothesis}
\newtheorem {thm}[hypo]{Theorem}   
\newtheorem {prop}[hypo]{Proposition}

\newtheorem {cor}[hypo]{Corollary}
\newtheorem {lem}[hypo]{Lemma}

\theoremstyle{remark}\newtheorem{rmk}[hypo]{Remark}

\begin{document}

\title[Symmetric cubic surfaces and  $\text{\rm G}_2$ character varieties]{Symmetric cubic surfaces and $\bf G_2$ character varieties}
\author{Philip Boalch and Robert Paluba}

\begin{abstract}
We will consider a two dimensional ``symmetric'' subfamily
of the four dimensional family of Fricke cubic surfaces.
The main result is that such symmetric cubic surfaces arise as character varieties for the exceptional group of type $G_2$.
Further, this symmetric family will be related to the fixed points of the triality automorphism of $\Spin(8)$, and an example involving the finite simple group of order $6048$ inside $G_2$ will be considered.
\end{abstract}

\maketitle

\section{Introduction}

The Fricke family of cubic surfaces is the family
\vspace{.3cm}
\beq\label{eq: fcs}
x\,y\,z+x^2+y^2+z^2 + b_1\,x+b_2\,y+b_3\,z + c =0
\eeq

\noindent
of affine cubic surfaces parameterised by constants
$(b_1,b_2,b_3,c)\in \IC^4$.
Note that the moduli space of (projective) cubic surfaces is four dimensional and a generic member will have an affine piece of this form, so this family includes an open subset of all cubic surfaces.
We will say a Fricke surface is {\em symmetric} if $b_1=b_2=b_3$.

The full family \eqref{eq: fcs} is known to be a semiuniversal deformation of a $D_4$ singularity (which occurs at the symmetric surface $b_i=-8,c=28$). 
Many other examples of cubic surfaces are isomorphic to symmetric Fricke cubics:
the Markov cubic surface ($b_i=c=0$), Cayley's nodal cubic surface  ($b_i=0, c=-4$), Clebsch's diagonal cubic surface ($b_i=0,c=-20$)
and the Klein cubic surface ($b_i=-1,c=0$), see \S\ref{sn: klein} below.

The Fricke surfaces are interesting since they are some of the simplest nontrivial examples of complex character varieties (and as such they are amongst the simplest examples of complete \hk 
four-manifolds for which we do not know how to construct the metric by finite dimensional means, cf. \cite{ihptalk} \S3.2).
Namely if $\Si=\IP^1\setminus\{a_1,a_2,a_3,a_4\}$ is a four
punctured sphere, then the moduli space
$$\MB(\Si,\SL_2(\IC)) = \Hom(\pi_1(\Si),\SL_2(\IC))/\SL_2(\IC)$$
of $\SL_2(\IC)$ representations of the fundamental group of $\Si$ 
is a (complex) six-dimensional algebraic Poisson variety and its symplectic leaves are Fricke cubic surfaces.
Indeed, choosing generators of the fundamental group leads to the identification of the character variety 
$$\MB(\Si,\SL_2(\IC)) \cong \SL_2(\IC)^3/\SL_2(\IC)$$
with the quotient of three copies of $\SL_2(\IC)$ by diagonal conjugation. (All our quotients are affine geometric invariant theory quotients, taking the variety associated to the ring of invariant functions.) 
In this case the ring of invariant functions is generated by the seven functions
$$ x= \tr(M_2M_3),\quad y= \tr(M_1M_3),\quad z=\tr(M_1M_2),$$
\beq\label{eq: fricke from sl2}
 m_1=\tr(M_1),\quad m_2=\tr(M_2),\quad m_3=\tr(M_3), \quad
m_4=\tr(M_1M_2M_3)
\eeq
where $M_i\in \SL_2(\IC)$.
These generators satisfy just one relation, given by 
equation \eqref{eq: fcs}, with:
$$ 
b_1 = -(m_1m_4 + m_2m_3), \ 
b_2 = -(m_2m_4 + m_1m_3), \ 
b_3 = -(m_3m_4 + m_1m_2), \ 
$$
\beq \label{eq: 2x2 formulae}
c=\Prod m_i -4+ \sum m_i^2.
\eeq
This relation amongst the generators is known as the Fricke relation\footnote{apparently (\cite{Magnus}) it  was discovered by Vogt (\cite{vogt} eq. (11)) in 1889, and repeatedly rediscovered by many others, including Fricke (\cite{FrickeKleinI} p.366).}.
The symplectic leaves are obtained by fixing the conjugacy classes of the monodromy around the four punctures and in general this amounts to fixing the values of the four invariants $m_1,m_2,m_3,m_4$, and thus the symplectic leaves are Fricke cubics.
Note that from this point of view the $D_4$ singularity occurs at the trivial representation of $\pi_1(\Si)$.

The aim of this article is to consider some simple examples of character varieties for the exceptional simple group $G_2(\IC)$ of dimension $14$.
Our main result (Corollary \ref{cor: main cor}) may be summarised as:

\begin{thm}\label{thm: 1}
There is a two parameter family of  character varieties for the exceptional group $G_2(\IC)$
which are isomorphic to smooth symmetric Fricke cubic surfaces, and thus to character varieties for the group $\SL_2(\IC)$.
\end{thm}

Note that the character varieties for complex reductive groups (and the naturally diffeomorphic Higgs bundle moduli spaces) are crucial to the geometric version of the Langlands program \cite{BD96}, and, in the case of compact curves $\Si$, the geometric Langlands story for $G_2$ shows significant qualitative differences to the case of $\SL_n$, involving a nontrivial involution of the  Hitchin base (see e.g. \cite{Hit-G2}). Thus it is surprising that there are $G_2$ character varieties which are in fact isomorphic to $\SL_2$ character varieties\footnote{We expect there to be analogous isomorphisms also in the Dolbeault/Higgs and De\! Rham algebraic structures corresponding to the Betti version considered here.}.  
As far as we know this is the first example of an isomorphism between nontrivial (symplectic) Betti moduli spaces involving an exceptional group.
Some further motivation is described at the end of section \ref{sn: charvars}.

Philosophically we would like to separate (or distance) the choice of the group from the ``abstract''  moduli space, as in the theory of Lie groups, where it is useful to consider the abstract group independently of a given representation, or embedding in another group.
In this language our result says that an abstract (nonabelian Hodge) moduli space has a ``$G_2$ representation'' as well as an $\SL_2$ representation/realisation.

In the later sections of the paper we will also consider the following topics. 
In \S\ref{sn: braid} the natural braid group actions on the spaces (coming from the mapping class group of the curve $\Si$) will be made explicit and we will show that the isomorphisms of Theorem \ref{thm: 1} are braid group equivariant. 
In \S\ref{sn: triality} we will recall that the $\IC^4$ parameter space of Fricke cubics is naturally related to the 
Cartan subalgebra of $\Spin_8(\IC)$ (i.e. the simply connected group of type $D_4$) and show that the subspace $\IC^2\subset \IC^4$ of symmetric Fricke cubics corresponds to the 
inclusion of Cartan subalgebras
$$ \lt_{G_2} \subset \lt_{\Spin(8)}$$
coming from the inclusion $G_2(\IC)\subset \Spin_8(\IC)$ identifying 
$G_2(\IC)$ as the fixed point subgroup of the triality automorphism of 
$\Spin_8(\IC)$.

\begin{figure}[h]
	\centering
	\input{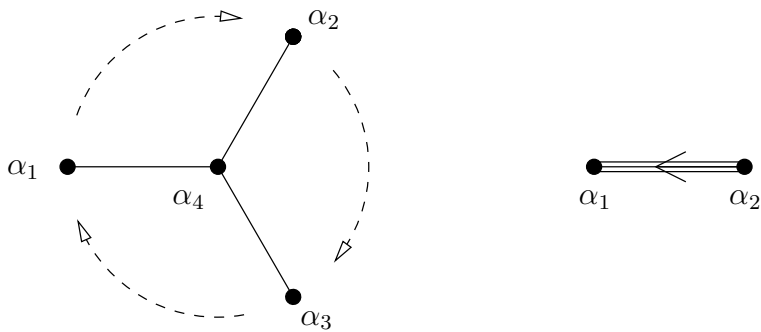}
	\caption{Triality automorphism of $D_4$ and the resulting $G_2$ Dynkin diagram}\label{fig: triality}
\end{figure}

Finally in section \S\ref{sn: klein} we will revisit some of the finite braid group orbits on cubic surfaces found in \cite{k2p-short,octa}.
In particular we will consider the Klein cubic surface (the unique cubic surface containing a braid group orbit of size $7$)
and show that:

\begin{thm}
If the Klein cubic surface $K$ is realised as a $G_2$ character variety (via Theorem \ref{thm: 1}) then the braid orbit of size $7$ in $K$ corresponds to some triples of generators of the finite simple group $G_2(\IF_2)'\subset G_2(\IC)$ of order $6048$.
One such triple of generators is uniquely determined by the three lines passing through a single point in the Fano plane $\IP^2(\IF_2)$. 
\end{thm}

\section{Tame character varieties}\label{sn: charvars}

Let $\Si$ be a smooth complex algebraic curve and let $G$ be a connected complex reductive group.
Then we may consider the space $$\Hom(\pi_1(\Si,p),G)$$ 
of representations
of the fundamental group of $\Si$ into the group $G$ (where $p\in \Si$ is a base point).
This is an affine variety with an action of $G$ given by conjugating representations. 
The character variety
$$\MB(\Si,G) = \Hom(\pi_1(\Si),G)/G$$
is defined to be the resulting affine geometric invariant theory quotient (the variety associated to the ring of $G$ invariant functions on 
$\Hom(\pi_1(\Si,p),G)$).
It is independent of the choice of basepoint so we supress $p$ from the notation.
Set-theoretically the points of $\MB(\Si,G)$ correspond bijectively to the 
{\em closed} $G$-orbits in $\Hom(\pi_1(\Si,p),G)$.

It is known that $\MB(\Si,G)$ has a natural algebraic Poisson structure \cite{AB83, Aud95long}.
The symplectic leaves of $\MB(\Si,G)$ are obtained as follows.
Suppose $\Si = \bar\Si\setminus \{a_1,\ldots a_m\}$ is obtained by removing $m$ points from a smooth compact curve $\bar \Si$.
Choose a conjugacy class $$\cC_i\subset G$$ 
for $i=1\ldots,m$ and let 
$\bcC$ denote this $m$-tuple of conjugacy classes.
Consider the subvariety
$$\MB(\Si,G,\bcC)  \subset \MB(\Si,G)$$
consisting of the representations taking a simple loop around $a_i$ into  
the class $\cC_i$ for each $i$.
Then the symplectic leaves of $\MB(\Si,G)$ are the connected components of 
the subvarieties $\MB(\Si,G,\bcC)$. We will also
refer to  these symplectic varieties $\MB(\Si,G,\bcC)$ (and their connected components)  as character varieties.
Note that in general $\MB(\Si,G,\bcC)$ is not an affine variety (i.e. it is not a closed subvariety of   $\MB(\Si,G)$), although it will be if  all the conjugacy classes $\cC_i\subset G$ are semisimple (since this implies that each $\cC_i$ is itself an affine variety), and this will be  the case in the examples we  will focus on here.

\subsection{Ansatz}

Consider the following class of examples of character varieties.
Let $G$ be as above and 
choose $n$ distinct complex numbers $a_1,\ldots,a_n$.
Let
$$\Si = \IC \setminus \{a_1,\ldots,a_n\} =  \IP^1(\IC) \setminus \{a_1,\ldots,a_n,\infty\}$$
be the $n$-punctured affine line (i.e. an $n+1$-punctured Riemann sphere),
and suppose conjugacy classes $\cC_1,\ldots,\cC_n,\cC_\infty\subset G$ are chosen so that:

a) $\cC_1,\ldots,\cC_n\subset G$ are semisimple conjugacy classes of minimal possible (positive) dimension, and

b) $\cC_\infty\subset G$ is a regular semisimple conjugacy class.

For example:

1) If $G=\SL_2(\IC)$ then these conditions just say that all the conjugacy classes are semisimple of dimension two,

2) If $G=\GL_n(\IC)$ then the minimal dimensional semisimple conjugacy classes are those with an eigenvalue of multiplicity $n-1$; after multiplying by an invertible scalar to make the corresponding eigenvalue $1$,
any such class contains complex reflections (i.e. linear automorphisms of the form ``one plus rank one''), and the monodromy group will be generated by $n$ complex reflections (note that the term  ``complex reflection group'' is often used to refer to such a group which is also {\em finite}),

3) If $G=G_2(\IC)$ then this means all the classes $\cC_1,\ldots,\cC_n$ are equal to the unique semisimple orbit $\cC\subset G$ of dimension $6$, and that $\cC_\infty$ is one of the twelve dimensional semisimple conjugacy classes. 

Suppose we now (and for the rest of the article) specialise to the case $n=3$, so that $\Si$ is a four punctured sphere.
This is the simplest genus zero case such that $\Si$ has an interesting  mapping class group.
Then it turns out that the above ansatz yields  two-dimensional character varieties in all of the above cases.

\begin{lem}\label{lem: dim 2 spaces}
Suppose $n=3$ so that $\Si$ is a four punctured sphere.
Then in each case $1),2),3)$ listed above the corresponding character variety $\MB(\Si,G,\bcC)$ 
is of complex dimension two, provided $\cC_\infty$ is sufficiently generic. 
\end{lem}
\pf
We will just sketch the idea: assuming $\cC_\infty$ is generic the action of $G$ on 
$$\{ (g_1,g_2,g_3,g_\infty)\st g_i\in \cC_i,\  g_1g_2g_3g_\infty=1\}$$
will have stabiliser of dimension $\dim(Z)$ where $Z\subset G$ is the 
centre of $G$.
The expected dimension of $\MB(\Si,G,\bcC)$ 
is then $\sum\dim(\cC_i) - 2\dim(G/Z)$.

For 1) each orbit is dimension two so it comes down to the sum
$$ 4\times 2 - 2\times 3 =2.$$
The resulting surfaces are the Fricke surfaces.

For 2) if $G=\GL_3(\IC)$ the first three orbits are of dimension four and the generic orbit is of dimension six, so the dimension of the character variety is
$$ 3\times 4 + 6 - 2\times (9-1)  =2.$$
This case was first analysed in \cite{pecr-nopham, k2p-short} and the resulting character varieties were explicitly related to the full four parameter family of Fricke surfaces (using the Fourier--Laplace transform).
In brief, the analogue of the Fricke relation that arises is given in \cite{k2p-short} (16) and this relation is related to the Fricke relation in  \cite{k2p-short} Theorem 1, using an explicit algebraic map that is derived from the Fourier--Laplace transform.

For 3), the group $G=G_2(\IC)$ has dimension $14$ (with trivial centre), the first three orbits are of dimension six and the generic orbit is of dimension twelve, so the dimension of the character variety is
$$ 3\times 6 + 12 - 2\times 14  =2$$
so is again a complex surface.
Note that since $G_2$ has rank two, there is a two-parameter family of choices for the orbit $\cC_\infty$ and so we expect to obtain a two parameter family of surfaces in this way.
\epf

\noindent
This yields our main question: what are the complex surfaces arising in the $G_2$ case?

Some further motivation for this question is as follows:

1) In \cite{ihptalk} \S3.2 there is a conjectural classification of the \hk manifolds of real dimension four that arise in nonabelian Hodge theory.
If this conjecture is true we should be able to locate the $G_2$ character varieties of Lemma \ref{lem: dim 2 spaces} on the list of \cite{ihptalk}. We do this here. (These complete \hk manifolds are noncompact analogues of K3 surfaces.)

2) The article \cite{logahoric} introduced the notion of parahoric bundles 
(by defining the notion of weight for parahoric torsors), the notion of logarithmic connection on a parahoric bundle (aka ``logahoric connections''), and established a precise Riemann--Hilbert correspondence for them. 
At first glance these seem to be quite exotic objects, and
the $G_2$ example considered here is one of the simplest contexts where this Riemann--Hilbert correspondence is necessary.
Our main result shows that the corresponding character varieties are in fact not at all exotic. 
The moduli spaces of such connections will be considered elsewhere.

\section{Octonions and $G_2$}

The compact exceptional simple Lie group $G_2$ arises as the group of automorphisms of the octonions. 
In this section we will describe the corresponding complex simple group $G_2(\IC)$, as the group of automorphisms of the complex octonions (or complex Cayley algebra)
and describe its unique semisimple conjugacy class $\cC\subset G_2(\IC)$ of dimension six.

This conjugacy class is somewhat remarkable since there is not a six dimensional semisimple adjoint orbit $\cO\subset \g_2=\Lie(G_2(\IC))$.

\begin{figure}[h]
	\centering
	\input{fano.plane.pstex_t}
	\caption{Points and lines in the Fano plane $\IP^2(\IF_2)$}\label{Fano Plane}
\end{figure}

\subsection{Complex octonions}

Let $\IO=\IO(\IC)$ denote the 8 dimensional non-associative algebra
with $\IC$-basis the symbols
$$1,e_1,e_2,e_3,e_4,e_5,e_6,e_7$$
and with multiplication table determined by
 the Fano plane $\IP^2(\IF_2)$ (and the fact that $1$ is central),
see Figure \ref{Fano Plane}.
Namely each triple $i,j,k$ of vertices of an oriented line in the Fano plane (for example $e_1,e_2,e_4$) forms a quaternionic triple:
$$i^2=j^2=k^2=ijk=-1.$$
Thus for example $e_1e_2=e_4$.
(Note, of course, that any two points lie on a unique line and so this determines the multiplication.) 
This multiplication table is symmetric under both of the two permutations
$$ e_n \mapsto e_{n+1}$$
and
$$ e_n \mapsto e_{2n}$$
of the indices (where all the indices are read modulo $7$, and we prefer to write $e_7$ rather than $e_0$).
This immediately implies that the triples of basis vectors with indices 
$$124, 235, 346, 457, 561, 672, 713$$
form quaternionic triples, so it is easy to remember the labelling on the Fano plane.
(Beware that some authors use a less symmetric multiplication table such that $e_1,e_2,e_3$ is a quaternionic triple; we follow 
\cite{baez-octonions, conway-smith}.)

Let $V\cong \IC^7$ denote the complex span of the 
$e_1,e_2,e_3,e_4,e_5,e_6,e_7$ and let
$$\tr : \IO\to \IC$$
denote the linear map with kernel $V$ taking $1\in\IO$ to $1\in \IC$.
Define an involution $q\mapsto \bar q$ of $\IO$ to be the $\IC$-linear map fixing $1$ and acting as $-1$ on $V$.
Then there is a nondegenerate symmetric $\IC$-bilinear form on $\IO$
$$\langle q_1,q_2 \rangle = \Tr(q_1\cdot\bar{q}_2)\in \IC.$$
Let $n(q) = \langle q,q \rangle\in \IC$ denote the corresponding quadratic form and we will say that $n(q)$ is the norm of $q$.
Note that if $v_1,v_2\in V$ then $\langle v_1,v_2 \rangle = -\Tr(v_1\cdot v_2)$.

The group $G_2(\IC)$ is the group of algebra automorphisms of $\IO$.
As such it acts trivially on $1$ and faithfully on $V$, preserving the quadratic form, so there is an embedding
$$G_2(\IC)\subset \SO(V,n) \cong \SO_7(\IC).$$ 
Henceforth we will write $G=G_2(\IC)$ and think of it in this 
seven dimensional representation.

\subsection{The six dimensional semisimple conjugacy class}

Let $T\subset G$ be a maximal torus, so that $T\cong (\IC^*)^2$.
Any semisimple conjugacy class in $G$ contains an element of $T$, and so to determine the dimension of the possible semisimple conjugacy classes it is sufficient to study the centralisers in $G$ of the elements of $T$.

Of course any element $t\in T$ is the exponential of some element
$X\in \lt = \Lie(T)$, but it is not always true that the centralisers
of $X$ and $t$ are the same and so care is needed (this is essentially the phenomenon of resonance in the theory of linear differential equations). This has been analysed in detail by Kac \cite{kac-1969} and it is possible to determine the centraliser of $t$ in terms of $X$ (see Serre \cite{serre-kac-coords}).
The result is that there is a unique semisimple conjugacy class of dimension six containing certain (special) order three elements of 
$G$. The centraliser of such elements is a copy of 
$\SL_3(\IC)\subset G$, and so we see immediately that 
$\dim \cC = \dim G/\SL_3(\IC)=14-8=6$.

In fact once we know (by the theorem of Borel--de\! Siebenthal or otherwise) that there is such an embedding of groups it is clear that the desired element of $G$ generates the centre of $\SL_3(\IC)$, and so is of order $3$.

We will skip the details of the above discussion since in terms of octonions we can be more explicit, as follows.
First we will recall Zorn's proposition.
Suppose $a\in \IO$ is a nonzero octonion and let 
$T_a : \IO\to \IO$ denote the linear map
$$T_a(q) = a.q.a^{-1}$$
given by conjugation by $a$.
(This is well-defined  since $\IO$ is associative on subalgebras generated by pairs of elements.)
However non-associativity implies such 
maps are not always automorphisms of $\IO$. 
Zorn characterised which maps $T_a$ are automorphisms:
\begin{prop}[Zorn, see \cite{conway-smith} p.98]
$T_a\in G$ if and only if $a^3\in \IC.1$
\end{prop}

Thus suppose we take an ``imaginary'' octonion $v\in V$ of norm 
$n(v)=3$. 
This means that $v$ is a square root of $-3$, i.e. $v^2=-3$.
Then we can consider the element
\beq\label{eq: av}
a(v) = \frac{1+v}{2}\in \IO.
\eeq
By construction $a(v)^3 = -1$, 
(cf. the fact that
 $(1+\sqrt{-3})/2 = \exp(\pi \sqrt{-1}/3)$).
Hence we have constructed an element 
$$T_{a(v)}\in G=G_2(\IC)$$
and it is clear that it is of order three in $G$.
The eigenvalues of $T_{a(v)}\in \Aut(V)$ are one 
(with multiplicity one) and the two nontrivial cube roots of one (each with multiplicity three). 

Let $\cO\subset V$ denote the set of elements of norm 3.
The group $G$ acts transitively on $\cO$, so $\cO$ is a single orbit of $G$ (in the representation $V$).

\begin{prop}\label{prop CeqO}
The map $\cO\to G=G_2(\IC)$ taking an element $v\in \cO$ to 
$T_{a(v)}\in G$ is a $G$-equivariant isomorphism of $\cO$
onto the six dimensional semisimple conjugacy class $\cC$ in $G$. 
\end{prop}
\pf
The $G$-equivariance is straightforward (where $G$ acts on itself by conjugation).
Thus $\cO$ is mapped onto a single conjugacy class.
To see the map is injective note that the eigenspace of $T_{a(v)}\in \Aut(V)$ with eigenvalue $1$ is one-dimensional and spanned by $v$, so it is sufficient to check that $T_{a(v)}\neq T_{a(-v)}$, but this is clear since they are inverse to each other (and of order three).
Note that $\cO$ is the quadric hypersurface $n(v)=3$ in $V$, so has dimension $6$
\epf

Note in particular that if $i,j,k$ is a quaternionic triple in $\IO$ then $i+j+k\in \cO$ and so we get an element 
$\frac{1}{2}(1+i+j+k)$
of $\cC$ for any line in the Fano plane. We will return to this in \S\ref{sn: klein}.

\section{Some invariant theory for $G_2$}

Our basic aim is to consider affine varieties obtained from the ring of $G=G_2(\IC)$ invariant functions on affine varieties of the form
\beq \label{eq: main quotient}
\{(g_1,g_2,g_3,g_\infty)\st g_1,g_2,g_3\in \cC, g_\infty\in \cC_\infty, g_1g_2g_3g_\infty=1\in G\}
\eeq
where $\cC\subset G$ is the six dimensional semisimple conjugacy class and $\cC_\infty$ is one of the twelve dimensional semisimple classes, and $G$ acts by diagonal conjugation.

Now a generic element $t$ of the maximal torus $T\subset G$ will be a member of such a conjugacy class $\cC_\infty$, and two such elements $t_1,t_2\in T$ are in the same class if and only if they are 
in an orbit of the action of the Weyl group $W =N_G(T)/T$ on $T$.

The Weyl group $W$ of $G$ is a dihedral group of order $12$ and its action on $T\cong (\IC^*)^2$ is well understood.
In brief there is a basis of $V$ such that $T$ is represented by diagonal matrices of the form
$$t=\diag(1,a_1,a_2,(a_1a_2)^{-1},a_1^{-1},a_2^{-1},a_1a_2)\in \GL(V)$$
for elements $a_1,a_2\in \IC^*$.
The action of $W$ on $T$ is generated by the two reflections
$$r_1(a_1,a_2) = (a_1^{-1},a_1a_2),\qquad 
  r_2(a_1,a_2) = (a_2,a_1)$$
fixing the hypertori $a_1=1$ and $a_1=a_2$ respectively. 

\begin{lem}[\cite{lorenz-mit} p.60]
The ring of $W$-invariant functions on $T$ is generated by the two functions 
$$\al = a_1+1/a_1+a_2+1/a_2 + a_1a_2 + 1/(a_1a_2)$$
$$\be = a_1/a_2+a_2/a_1+a_1^2a_2+a_1a_2^2 + 1/(a_1^2a_2) + 1/(a_1a_2^2).$$
\end{lem}

Note that if $t\in T$ is represented as a diagonal matrix as above then
$$
\al  = \tr_V(t) -1,\qquad
2\be =  \al^2-2\al-\tr_V(t^2)-5$$
so that specifying the values of $\al,\be$ is equivalent to specifying the values of $\tr_V(t)$ and $\tr_V(t^2)$. 
Of course the functions $\tr_V(t),\tr_V(t^2)$ are just the restrictions of the functions $\tr_V(g),\tr_V(g^2)$ defined on $G\subset \GL(V)$, and so (in this way) we can just as well view $\al$ and $\be$ as invariant functions on $G$. (The resulting formulae are simpler if we work with $\be$ rather than $\tr_V(g^2)$.)

Thus we may rephrase our main question differently:
suppose we consider the affine variety
$$\IM := \cC^3/G$$
associated to the ring of $G$-invariant functions on $\cC^3$ (where $G$ acts by diagonal conjugation).
Then the affine varieties 
$$\{(g_1,g_2,g_3,g_\infty)\st g_1,g_2,g_3\in \cC, g_\infty\in \cC_\infty, g_1g_2g_3g_\infty=1\in G\}/G$$
(associated to the invariant functions on \eqref{eq: main quotient})
arise as fibres of the map
$$\pi: \IM \to \IC^2; \quad
[(g_1,g_2,g_3)]\mapsto (\al(g_1g_2g_3), \be(g_1g_2g_3)),$$
since fixing (sufficiently generic) values of the map $\pi$ will fix the conjugacy class of the product $g_1g_2g_3$, as required.
Thus as a first step we need to understand the affine variety $\IM = \cC^3/G$ and secondly we need to understand the map $\pi$ from $\IM$ to $\IC^2$.

\begin{prop}\label{prop: mc4}
The affine variety $\IM = \cC^3/G$ is isomorphic to an affine space of dimension four (more precisely the ring of $G$ invariant functions on $\cC^3$ is a polynomial algebra in four variables).
\end{prop}
\pf
Via Proposition \ref{prop CeqO} the affine variety $\cC^3$ is $G$-equivariantly isomorphic to $\cO^3$ where $\cO\subset V$ is the affine quadric of ``imaginary'' octonions of norm $3$. 

In other words $\cO^3$ is the subset of $(v_1,v_2,v_3)\in V^3$ such that  $n(v_i) = 3$ for $i=1,2,3$.
Thus we can first consider the quotient $V^3/G$.

Now G. Schwarz has classified all the representations of simple algebraic groups such that the ring of invariants is a polynomial algebra \cite{schwarz-regular}. 
Row 1 of Table 5a of \cite{schwarz-regular} (see also \cite{schwarz-g2bams}) says that 
the ring $\IC[V^3]^G$ is a polynomial algebra with 7 generators, and the generators can be taken to be the invariant functions:
\vspace{.3cm}
\begin{alignat}{1}
n(v_1),\quad &n(v_2),\quad n(v_3) \notag\\
p_1=\langle v_2,v_3\rangle,\quad 
p_2=&\langle v_1,v_3\rangle,\quad  
p_3=\langle v_1,v_2\rangle, \label{eq: pinvts}\\
p_4 = \langle &v_1,v_2.v_3\rangle.\notag
\end{alignat}

We are interested in the algebra obtained by fixing the values of the first three invariants to be $3$, which is thus a polynomial algebra generated by the remaining four functions $p_1,p_2,p_3,p_4$.
\epf

In other words the map 
$p=(p_1,p_2,p_3,p_4):\IM\to \IC^4$, with components the four functions $p_i$, identifies $\IM$ with $\IC^4$.
Thus the question now is to understand the map $\pi:\IM\to \IC^2$ in terms of the functions $p_i$.
The key formulae are the following, which may be verified symbolically.

\begin{thm}\label{thm: albe polys}
Suppose $v_1,v_2,v_3\in \cO\subset V\subset \IO$ are octonions in $V$ with norm $3$, and $g_i = T_{a(v_i)}\in \cC\subset G\  (i=1,2,3)$ are the corresponding elements in the six dimensional semisimple conjugacy class in $G=G_2(\IC)$.
Then the invariant functions
$$\al=\al(g_1g_2g_3),\qquad \be=\be(g_1g_2g_3)$$
of the product $g_1g_2g_3\in G$
are related to the invariants
$$
p_1=\langle v_2,v_3\rangle,\quad 
p_2=\langle v_1,v_3\rangle,\quad  
p_3=\langle v_1,v_2\rangle,\quad
p_4 = \langle v_1,v_2.v_3\rangle$$ 
of the octonions $v_1,v_2,v_3$ by the following formulae:
$$8\,\al = 
{ p_4}\,{ s_1}-{{ s_1}}^{2}+3\,{ p_4}+3\,{ s_1}+3\,{ s_2}+{ s_3}-6,$$
\begin{align*}    %
64\,\be = 
-{{ p_4}}^{3}&+3\,{{ p_4}}^{2}{ s_1}-3\,{ p_4}\,{{ s_1}}^{2}-
7\,{{ s_1}}^{3}+9\,{{ p_4}}^{2}-12\,{ p_4}\,{ s_1}+18\,{ p_4
}\,{ s_2}+6\,{ p_4}\,{ s_3}\\
&+39\,{{ s_1}}^{2}+18\,{ s_1}\,{
 s_2}+6\,{ s_1}\,{ s_3}-9\,{ p_4}-9\,{ s_1}-90\,{ s_2}-30
\,{ s_3}-183,
\end{align*}

\noindent
where
$s_1=p_1+p_2+p_3,\ 
s_2=p_1p_2+p_2p_3+p_3p_1,\ 
s_3=p_1p_2p_3$.

\end{thm}

Now we will explain how the symmetric Fricke cubics arise.
Fix constants $k_\al,k_\be\in \IC$ and consider the subvariety
of $\IM\cong \IC^4$ cut out by the equations
$$\al=k_\al,\qquad \be=k_\be$$
where $\al,\be$ are viewed as functions on $\IM$ via the above formulae.

First suppose we change coordinates on $\IM$ by replacing $p_4$ by the function $b$ defined so that
\beq\label{eq: bdefn}
p_1+p_2+p_3+p_4 = 4b+5.
\eeq
Then, in the presence of the relation $\al=k_\al$, the equation
$\be=k_\be$ simplifies to the equation:
\beq\label{eq: bcubic}
{b}^{3}+6\,{b}^{2}-3\,({k_\al}-1)\,b+{k_\be}+2=0.
\eeq
This is a cubic equation for $b$ and so specifying $k_\al,k_\be$ determines $b$ up to three choices.

Now we can reconsider the remaining equation $\al=k_\al$, which is easily seen to be a symmetric Fricke cubic:
expanding the symmetric functions and making the substitutions
\beq\label{eq: p to xyz}
p_1 = 1-2x,\quad p_2=1-2y,\quad p_3=1-2z
\eeq
converts the equation $\al=k_\al$ into the Fricke cubic
\beq\label{eq: symfcs}
x\,y\,z + x^2+y^2+z^2 + b\,(x+ y+ z) + c =0
\eeq
where 
\beq\label{eq: c formula}
c=k_\al-2-3\,b.
\eeq

The situation may be  summarised in the following commutative diagram:

\begin{figure}[h]
	\centering
	\input{square.pstex_t}
	\caption{Main diagram}\label{fig: main diag}
\end{figure}

Here, in the top left corner,  $\IM=\cC^3/G_2(\IC)$ which we identified with $\IC^4$ using the functions $p_1,p_2,p_3,p_4$ in Proposition \ref{prop: mc4}.
The map $\pi$ on the left is given by the invariant functions $\al,\be$, expressed as explicit functions on $\IM$ as in Theorem \ref{thm: albe polys}. The fibre of the map $\pi$ over a generic point $(k_\al,k_\be)\in\IC^2$ is the $G_2(\IC)$ character variety
$\MB(\Si,G_2(\IC),\bcC)$ where $\bcC=(\cC,\cC,\cC,\cC_\infty)$ 
with $\cC_\infty$ the twelve dimensional conjugacy class in $G_2(\IC)$  with eigenvalues having invariants $\al=k_\al,\be=k_\be$.

On the right we consider a universal family of symmetric Fricke cubics: we take a copy of $\IC^4$ with coordinates $x,y,z,b$ and consider the map $\wh \pi$ to $\IC^2$ given by $(b,c)$, where $c$ is viewed as a function of $x,y,z,b$ via the Fricke equation \eqref{eq: symfcs}.
The fibre of $\wh \pi$ over a point $(b,c)\in \IC^2$ is the symmetric Fricke cubic \eqref{eq: symfcs} with these values of $b,c$.

The isomorphism $\varphi$ along the top is defined by the equations
\eqref{eq: bdefn},\eqref{eq: p to xyz} expressing $x,y,z,b$ in
terms of $p_1,p_2,p_3,p_4$.

The map $\pr$ along the bottom is given by the equations 
\eqref{eq: bcubic},\eqref{eq: c formula}
expressing the values $k_\al,k_\be$ of $\al,\be$ in terms of $b,c$:
\beq\label{eq: pr formula}
\al= c+2+3\,b,\qquad
\be=-{b}^{3}+3\,{b}^{2}+3\,bc+3\,b-2.
\eeq
This is a finite surjective map, with generic fibres consisting of three points. In fact by examining the discriminants we can be more precise.

\subsection{Discriminants}
First we can consider the discriminant $\cD_{\pr}\subset \IC^2_{\al,\be}$  
of the cubic polynomial 
\eqref{eq: bcubic}, characterising when the fibres of $\pr$ have less than $3$ elements. This discriminant locus (with $\al,\be $ replacing $k_\al,k_\be$ resp.) is: 
\beq\label{eq: bdisc}
4 \al^3-12 \al \be-\be^2-36\al-24\be-36=0.
\eeq

Secondly we can consider the discriminant 
$\cD_W\subset \IC^2_{\al,\be}$
of the quotient map 
$(\al,\be):T\to T/W\cong \IC^2$. 
This discriminant is the image of the $6$ ``reflection hypertori'' in $T$ and is the subvariety of $\IC^2$ for which the fibres of this quotient map have less than $12=\#W$ points.
Explicitly  it is cut out by the square of the Weyl denominator function, which  (\cite{FulHar} p.413) is the following $W$-invariant function on $T$:
$$\bigl((a_1-a_2)(a_2-a_3)(a_3-a_1)\bigr)^2 \times 
(a_1a_2+a_2a_3 + a_3a_1- a_1-a_2-a_3)^2$$
where $a_3 = 1/(a_1a_2)$.   
The two factors here correspond to the long and the short roots of $G_2$.
If we now rewrite these two factors in terms of the basic invariants $\al,\be$ we find that $\cD_W=\cD_1\union\cD_2$ has two irreducible components corresponding to the two factors, and that the first component is equal to \eqref{eq: bdisc}, the discriminant locus  $\cD_{\pr}$
 of the map $\pr$. 
The other irreducible component $\cD_2$ of $\cD_W$ is given by
\beq\label{eq: d2}
{{\al}}^{2}-4\,{\be}-12=0.
\eeq
Since $\cD_{\pr}\subset \cD_W$ 
this implies that: {\em if $\cC_\infty$ has dimension twelve then the
corresponding values $k_\al,k_\be$ of the invariants $\al,\be$ are not on the discriminant \eqref{eq: bdisc}, and so the fibre $\pr^{-1}(k_\al,k_\be)$ has exactly three points}.

In turn we can relate this to the locus of {\em singular} symmetric Fricke cubics surfaces.

\subsection{Singular points}

Given $b,c\in \IC$ the corresponding symmetric Fricke surface 
\eqref{eq: symfcs} is singular if and only if the derivatives of the defining equation all vanish, i.e. if there is a simultaneous solution $(x,y,z)$ of \eqref{eq: symfcs} and the three equations
$$xy+2z+b=0,\qquad yz+2x+b=0,\qquad xz+2y+b=0.$$
Eliminating two variables and considering the resultant of the remaining two equations we find:

\begin{lem}
The symmetric Fricke surface \eqref{eq: symfcs} determined by $b,c\in \IC$ is singular if and only if
\beq\label{eq: sing locus}
\left( {b}^{2}-8\,b-4\,c-16 \right)  \left( 4\,{b}^{3}-3\,{b}^{2}-6\,
bc+{c}^{2}+4\,c \right)=0.
\eeq
\end{lem}

Now we can consider the image of this singular locus 
$\cD_{sing}\subset \IC^2_{b,c}$ under the map $\pr:\IC^2_{b,c}\to\IC^2_{\al,\be} $.

\begin{prop}\label{prop: mapsinglocus}
The singular locus $\cD_{sing}\subset \IC^2_{b,c}$ maps onto the discriminant locus 
$\cD_W\subset \IC^2_{b,c}$. 
Specifically the first (left-hand) irreducible 
component of $\cD_{sing}$
maps onto the first component $\cD_1=\cD_{\pr}$ of $\cD_W$ given in \eqref{eq: bdisc}, and the second component of $\cD_{sing}$
maps onto the second component $\cD_2$ of $\cD_W$, given in \eqref{eq: d2}.
\end{prop}

Note however that the inverse image $\pr^{-1}(\cD_1)$ has another irreducible component 
\beq\label{eq: dbl locus}
{b}^{2}+b-c-1=0
\eeq
(with multiplicity two), besides the first component of the singular locus \eqref{eq: sing locus}.

Consequently we see that:

\begin{cor}\label{cor: main cor}
Let $\cC\subset G_2(\IC)$ be the six dimensional semisimple conjugacy class.
Then for any regular semisimple conjugacy class  
$\cC_\infty\subset G_2(\IC)$ the character 
variety $$\MB(\Si,G_2(\IC),\bcC)$$ 
with $\bcC=(\cC,\cC,\cC,\cC_\infty)$, 
has three connected components, each of which is isomorphic to a smooth symmetric Fricke cubic surface.
\end{cor}
\pf
Such a character variety is the fibre of the map $\pi$ over a point 
$(k_\al,k_\be)\in \IC^2\setminus\cD_W$.
Since $\pi$ factors through $\pr$ (going around the square, i.e. via the dashed diagonal map), and $\cD_{\pr}\subset \cD_W$, such fibres consist of three fibres of $\wh \pi$ over the smooth locus $\IC_{b,c}^2\setminus \cD_{sing}$. These fibres are smooth symmetric Fricke cubic surfaces.
\epf

Note in particular this implies we can canonically associate two other smooth cubic surfaces to any smooth symmetric Fricke cubic surface with parameters $b,c$ not on the conic \eqref{eq: dbl locus}, namely  the other two components of  $\pi^{-1}(\pr(b,c))\subset \IM$.

Note also that in \cite{k2p-short} the map relating the $\GL_3(\IC)$ character varieties to the Fricke cubics (as in part 2 of Lemma \ref{lem: dim 2 spaces} above) was derived from earlier work on the Fourier--Laplace transform, 
whereas here we have identified the character varieties directly;
we do not (yet) understand if there is a $G_2$ analogue/extension of Fourier--Laplace.

In the next section we will consider the natural braid group actions and show that the isomorphism $\varphi$ is braid group equivariant.

\begin{rmk}
Suppose instead we replace  $\cC_\infty$ by the closure of the regular unipotent conjugacy class in the definition of $\MB(\Si,G_2(\IC),\bcC)$.
The resulting variety is the fibre of $\pi$ over the point $\al=\be=6$, and we readily see it has two components: the Fricke cubic $(b,c)=(-8,28)$ with the $D_4$ singularity, and (with multiplicity two) the  Fricke cubic with 
$(b,c)=(1,1)$. This surface is also singular; in fact these
 parameters lie at 
the cusp of the cuspidal cubic on the right in the singular locus \eqref{eq: sing locus}. 
\end{rmk}

\begin{rmk}
The first irreducible component of  the singular locus 
\eqref{eq: sing locus} could be called the ``very symmetric locus'' of Fricke cubics, since in the original Fricke--Vogt $\SL_2(\IC)$ 
picture it corresponds to the case 
where {\em all} the local monodromies are conjugate: $m_1=m_2=m_3=m_4$
so that, by the formulae \eqref{eq: 2x2 formulae}: 
$$ b=-2m^2,\qquad c = m^4-4+4m^2 $$
and thus ${b}^{2}-8b=4c+16$.
Although all these surfaces are singular this case has many applications, for example to anti-self-dual four-manifolds (\cite{Hit-tei} Theorem 3).
\end{rmk}

\begin{rmk}
This very symmetric  locus is closely related to the one-parameter family $b=0$  (which does admit smooth members, so is related to $G_2$ via Corollary \ref{cor: main cor}). 
This family is the family of $\SL_2(\IC)$-character varieties for the once-punctured torus (\cite{goldman-top} p.584, or \cite{vogt, Magnus}).
The basic statement relating this case to the very symmetric case is that there is a (degree four) ramified covering map between the two families, as follows (cf. \cite{k2p-short} Remark 14).
Suppose $b\in \IC$ and $c=b^2/4-2b-4$ so that the polynomial
$$ f=xyz+x^2+y^2+z^2 +b (x+y+z) + c$$
defines a very symmetric Fricke cubic.
If we define $d=-4-b/2$ and consider the polynomial
$$ g = XYZ+X^2+Y^2+Z^2 +d$$
then, if the variables are related by
$x=2-X^2, y=2-Y^2, z=2-Z^2$, the relation 
$$ f(x,y,z)  = g(X,Y,Z)g(-X,Y,Z)$$
holds, so the surface $g=0$ maps to the surface $f=0$, and the generic fibre contains $4$ points. 
Note that if $b=0$ then $c=d=-4$ so this is an endomorphism of the Cayley cubic surface. 
\end{rmk}

\section{Braid group actions}\label{sn: braid}

Since the symmetric Fricke cubics are character varieties they admit braid group actions, coming from the mapping class group of the four-punctured sphere $\Si$. In explicit terms (cf. \cite{icosa-short} \S4) this can be expressed in terms of changing the choice of loops generating the fundamental group of $\Si$, whence it becomes the classical Hurwitz braid group action:

\begin{lem}
Let $G$ be a group. Then there is an action of the three string Artin braid group $B_3$ on $G^3$, generated by the two operations $\be_1,\be_2$ where:
\begin{alignat}{1}
\be_1(g_1,g_2,g_3) &= (g_2,g_2^{-1}g_1g_2, g_3) \label{eq: be1}\\
\be_2(g_1,g_2,g_3) &= (g_1,g_3,g_3^{-1}g_2g_3). \label{eq: be2}
\end{alignat}
\end{lem}

Taking $G$ to be a complex reductive group, 
this describes the $B_3$ action on the character variety
$\MB(\Si,G)$, and given four conjugacy classes $\bcC\subset G^4$ this action restricts to the symplectic leaves
$$\MB(\Si,G,\bcC) \cong 
\{(g_1,g_2,g_3,g_\infty)\st g_i\in \cC_i, g_1g_2g_3g_\infty=1\}/G$$
provided the first three conjugacy classes are equal.

In the case $G=\SL_2(\IC)$ it is easy to compute the resulting action on the $G$-invariant functions on $\Hom(\pi_1(\Si),G)$ and in turn on the family of Fricke surfaces \cite{Iwas-modular,k2p-short}. 
In the symmetric case, for fixed constants $b,c\in \IC$ the formula is as follows:

\beq\label{eq: frick braiding}
\be_1(x,y,z) = (x,-b - z-xy,y),\qquad
\be_2(x,y,z) = (z,y,-b-x-yz).
\eeq 

The main aim of this section is to compute directly the action in the case
$G=G_2(\IC)$ we have been studying.
Let $\cC\subset G_2(\IC)$ denote the six dimensional semisimple conjugacy class, and let $\cO\subset V$ denote the orbit of elements of norm $3$.
\begin{prop}\label{prop: g2 braiding}
1) The braid group action \eqref{eq: be1},\eqref{eq: be2} on triples
$(g_1,g_2,g_3)\in \cC^3$ of elements of $\cC$ corresponds to the action
$$\be_1(v_1,v_2,v_3) = (v_2, \bar{w}_2\cdot v_1 \cdot {w}_2,v_3),$$
$$\be_2(v_1,v_2,v_3) = (v_1, v_3, \bar{w}_3\cdot v_2 \cdot{w}_3)$$
on triples of elements $(v_1,v_2,v_3)\in \cO^3\subset V^3\subset \IO^3$, via the isomorphism $\cC^3\cong\cO^3$ of Proposition \ref{prop CeqO}, where $w_i = (1+v_i)/2\in \IO$ (so that $g_i=T_{w_i}$).

2) The resulting $B_3$ action on $\IM=\cC^3/G\cong \IC^4$ is given by the formulae:
$$\be_1(p_1,p_2,p_3,p_4) = 
((p_4+p_1p_3-p_2)/2, p_1, p_3, (p_4+3p_2-p_1p_3)/2),$$
$$\be_2(p_1,p_2,p_3,p_4) = 
(p_1,(p_4+p_1p_2-p_3)/2, p_2, (p_4+3p_3-p_1p_2)/2)$$
in terms of the invariant functions $p_1,p_2,p_3,p_4$ defined in \eqref{eq: pinvts}. 
\end{prop}
\pf
1) is straightforward (taking care due to the non-associativity of $\IO$), and can be verified symbolically.

2) is now a nice exercise in expanding $G$-invariant functions on $\cO^3$ in terms of the basic invariants $p_i$. To illustrate this we will show $\langle v_3,(v_2v_1)v_2\rangle = 3p_2 -2p_1p_3$.
If $q\in \IO$ we will write $q=\Tr(q) +\Im(q)$ with 
$\Im(q)\in V$ and $\Tr(q)\in \IC.1$, and note that
$\langle v_1,v_2v_3\rangle$ is a skew-symmetric $3$-form if $v_1,v_2,v_3\in V$.
Then we have
\begin{alignat*}{1}
\langle v_3,(v_2v_1)v_2\rangle &= 
\langle v_3,\Tr(v_2v_1)v_2\rangle  +
\langle v_3,\Im(v_2v_1)v_2\rangle\\
&= -p_1p_3 - \langle v_3,v_2\Im(v_2v_1)\rangle \\
&= -p_1p_3 - \langle v_3,v_2(v_2v_1)\rangle + \langle v_3,v_2\tr(v_2v_1)\rangle\\
&= -p_1p_3 + 3p_2 -p_1p_3. 
\end{alignat*}
The stated formulae can all be derived in this way.
\epf

Observe that the sum $p_1+p_2+p_3+p_4$ is preserved by this braid group action. This corresponds to the fact that $b$ is preserved by the action on Fricke surfaces, and in fact there is a complete correspondence:

\begin{cor}\label{cor: breq}
The isomorphism  $\varphi$ in Figure \ref{fig: main diag}
identifying  $\IM=\cC^3/G$ with the universal family of symmetric Fricke cubic surfaces, 
is  braid group equivariant.
\end{cor}
\pf
Given the explicit formulae (in \eqref{eq: frick braiding} and Proposition \ref{prop: g2 braiding}) for the braid group actions, 
this follows directly from the formulae
\eqref{eq: bdefn},\eqref{eq: p to xyz} defining $\varphi$.
\epf

\section{Triality}\label{sn: triality}

The fact that the full family of Fricke cubics is a semiuniversal deformation of a $D_4$ singularity is not just a coincidence, and has deep moduli theoretic meaning\footnote{here we mean cubic surfaces as moduli spaces of local systems---see \cite{naruki-cr-looij} for more on the direct relation between $D_4$ and the moduli of cubic surfaces themselves.}. 
In brief the whole family of these Fricke moduli spaces (in fact the whole ``\hk nonabelian Hodge structure'') is naturally attached to the $D_4$ root system (or more precisely to the affine $D_4$ root system, but that is determined by the finite root system).
Note that the $D_4$ Dynkin diagram is the most symmetric (finite) Dynkin diagram since it has an automorphism of order three, the triality automorphism,  indicated in Figure \ref{fig: triality}. 

There are various ways to see the appearance of $D_4$ from the moduli space of rank $2$ representations of the fundamental group of the four-punctured sphere $\Si$:

1) Translating to our notation, Okamoto \cite{OkaPVI} found that if we label the eigenvalues of $M_i\in \SL_2(\IC)$ as $\exp(\pm \pi\sqrt{-1} \th_i)$ where $\th=(\th_1,\th_2,\th_3,\th_4)\in \IC^4$,
then this copy of $\IC^4$ should be viewed as the Cartan subalgebra of type $D_4$:   
he showed there is a natural action of the affine $D_4$ Weyl group on this $\IC^4$, and it lifts to automorphisms of the corresponding moduli spaces of rank two logarithmic connections on $\Si$ fibred over $\IC^4$ (at least off of the affine root hyperplanes---cf. \cite{AL-p6ims}).
Further 
Okamoto showed one can add in the automorphisms of the affine $D_4$ Dynkin diagram, to obtain an action of $\Sym_4\sdp W_{\text{aff}}(D_4)$ (which is isomorphic to the affine $F_4$ Weyl group).

2) The moduli spaces of rank two logarithmic connections on $\Si$ mentioned in 1) above have simple open pieces $\cM^*$ (where the underlying vector bundle on $\IP^1$ is holomorphically trivial)
and these open pieces are isomorphic to the $D_4$ asymptotically locally Euclidean \hk four-manifolds constructed by Kronheimer \cite{Kron.ale}. 
In fact Kronheimer constructs these spaces as a finite dimensional \hk quotient starting with a vector space of linear maps in both directions along the edges of the affine $D_4$ Dynkin graph (an early example of a ``quiver variety'').

3) The space of representations of the fundamental group of $\Si$ is closely related to the perverse sheaves on $\IP^1$ with singularities at the marked points, and such perverse sheaves have a well-known quiver description, again in terms of linear maps in both directions along the edges of an affine $D_4$ Dynkin graph.
This leads to the statement that the Fricke cubic surfaces are affine $D_4$ ``multiplicative quiver varieties'', in the sense of \cite{CB-Shaw, yamakawa-mpa}.

Note that each of these points of view has an extension to many other moduli spaces, often of higher dimensions (cf. \cite{rsode, slims-short}).

In any case, whichever is the reader's preferred viewpoint, 
it is natural to consider the action of the triality automorphism
on the space $\IC^4$ of parameters $\th$ (and the induced action on the Fricke coefficients $b_1,b_2,b_3,c$).
It is well-known that the root system of $G_2$ arises by `folding' the $D_4$ root system via the triality automorphism in this way, as indicated in Figure \ref{fig: triality}.
It turns out that the fixed locus is the space of symmetric Fricke cubics, so there is another (a priori different) link between symmetric Fricke cubics and $G_2$:

\begin{prop}
The action of the triality automorphism on $\th\in \IC^4$
permutes $\th_1,\th_2,\th_3$ cyclically and fixes $\th_4$.
The fixed locus $\th_1=\th_2=\th_3$ maps  onto the 
parameter space $b_1=b_2=b_3$  of the symmetric Fricke cubic surfaces.
\end{prop}
\pf
Note that the triality automorphism is only well defined upto conjugation by an inner automorphism, so we have some freedom. 
The key point is that, in Okamoto's setup, 
one should not use the standard (Bourbaki) convention for the $D_4$ root system, rather, as explained in \cite{icosa-short} Remark 5, it naturally appears as the set of {\em short} $F_4$ roots, i.e. as the set of $24$ norm one vectors:
\beq\label{eq: D4m}
\th\in D_4^-=\Bigl\{\pm \varepsilon_1,\quad\pm \varepsilon_2,\quad\pm \varepsilon_3,\quad\pm \varepsilon_4,\quad\frac{1}{2}(\pm \varepsilon_1\pm \varepsilon_2\pm \varepsilon_3\pm \varepsilon_4)\Bigr\} \subset \IC^4,
\eeq
where $\varepsilon_i$ are the standard (orthonormal) 
basis vectors of $\IC^4$ (i.e. as the group of unit  Hurwitz integral quaternions).
Then we just note that the vectors 
$$
\al_1 = \varepsilon_1,\quad
\al_2 = \varepsilon_2,\quad
\al_3 = \varepsilon_3,\quad
\al_4 = \frac{1}{2}(
\varepsilon_4-\varepsilon_1-\varepsilon_2-\varepsilon_3)$$
form a basis of simple roots, and the longest root is 
$\varepsilon_4 = \al_1+\al_2+\al_3+2\al_4$ 
and the mutual inner products of these simple roots are as indicated in the $D_4$ Dynkin diagram in Figure \ref{fig: triality}.
Thus the triality automorphism shown in Figure \ref{fig: triality}
acts as 
$$\tau(\th_1,\th_2,\th_3,\th_4) \ = \ (\th_3,\th_1,\th_2,\th_4)$$
and the fixed point locus is indeed $\th_1=\th_2=\th_3$.
Using the formulae \eqref{eq: fricke from sl2} it is clear that this maps to the locus of symmetric Fricke cubics.
\epf

\begin{rmk}
In fact the Fricke functions $b_1,b_2,b_3,c$ of $\th$ have a direct Lie-theoretic interpretation, as follows (this a minor modification of \cite{oblomkov-cubics} p.888-9, adjusting the weight, and thus coroot, lattice):
Let $\lt\cong\IC^4$ be the $D^-_4$ Cartan subalgebra.
Then the ``Fricke map'' 
\eqref{eq: fricke from sl2},\eqref{eq: 2x2 formulae}:
$$ \lt \to \IC^4; \qquad\th=\sum\th_i\varepsilon_i\mapsto (b_1,b_2,b_3,c)$$
is the quotient by the $D_4^-$ affine Weyl group (the semidirect product of the finite Weyl group $W$ and the coroot lattice $\Gamma_R$).
Indeed with these conventions (cf. also \cite{srops-short}) $\Gamma_R=\langle 2D_4^-\rangle\subset \lt$ and the quotient $\lt/\Gamma_R$ is a maximal torus $T_{\Spin(8)}$ of the simply connected group $\Spin_8(\IC)$ of type $D_4$.
In turn it is well known that the $W$-invariant functions on $T_{\Spin(8)}$ 
form a polynomial algebra generated by the fundamental weights  (\cite{FulHar} p.376), and these 
weights  $D_1,D_2,D^+,D^-$ are described in \cite{FulHar} (23.30): a straightforward computation then shows that
$$ D_1 =- b_1,\quad D_2=c,\quad D^+= -b_3,\quad D^{-} =-b_2$$
so it is clear that permuting the $b_i$ corresponds to triality, permuting the standard representation and the two half-spin representations of $\Spin_8(\IC)$.
\end{rmk}

Note that Manin \cite{Manin-P6} \S1.6 considered ``Landin transforms'' in this context. They are defined on a distinguished two-dimensional subspace of the full space $\IC^4$ of parameters: this looks to be different (and inequivalent) to the symmetric subspace we are considering (it looks to be the fixed locus arising from the involution of the affine $D_4$ Dynkin diagram swapping two pairs of feet, rather than the triality automorphism).

\section{The Klein Cubic Surface}\label{sn: klein}

The article \cite{k2p-short} found that there was a Fricke cubic surface containing a braid group orbit of size seven.
This arose by considering Klein's simple group of order 
$168=2^3\cdot3\cdot7$, 
the group of automorphisms of Klein's quartic curve (the modular curve $X(7)$):
$$X^3Y+ Y^3Z+Z^3X=0,$$
and then taking the corresponding complex reflection group in $\GL_3(\IC)$, of order $336$.
The braid group orbit of the conjugacy class of the standard triple of generating reflections of this complex reflection group has size $7$ and lives in a character variety of dimension two (as in part 2 of Lemma \ref{lem: dim 2 spaces} above).
Then using the Fourier--Laplace transform it was shown how to relate this $\GL_3(\IC)$ character variety to the usual $\SL_2(\IC)$ Fricke--Vogt story, i.e. to show that it is a Fricke cubic surface.
The resulting surface has 
$$m_1=m_2=m_3=2\cos(2\pi/7),\quad m_4=2\cos(4\pi/7)$$
(\cite{k2p-short} p.177) 
and the corresponding $\SL_2(\IC)$ monodromy group is a lift to $\SL_2(\IC)$ of the  (infinite) $2,3,7$ triangle group $\Delta_{237}\subset \PSL_2(\IC)$ (cf. \cite{octa} Appendix B)\footnote{in fact there are three inequivalent choices of seventh root of unity that one can make: for  two choices the projective monodromy group is a subgroup of $\PSU_2$ isomorphic to $\Delta_{237}$  and for the other choice one obtains the usual $\Delta_{237}\subset \PSL_2(\IR)$.}.
Thus from the formulae \eqref{eq: fricke from sl2}, and the fact that 
$4\cos(2\pi/7)\cos(4\pi/7)+4\cos(2\pi/7)^2=1$, we see that the 
{\em Klein cubic surface} is:
\beq\label{eq: klein cubic}
x\,y\,z+x^2+y^2+z^2=x+y+z
\eeq

\noindent
i.e. it is the ({smooth}) symmetric Fricke cubic with $b=-1,c=0$.
The braid orbit of size seven consists of the points: 
$$(x,y,z) = 
(0,0,0),\ 
(0,0,1),\ 
(0,1,0),\ 
(1,0,0),\ 
(1,1,0),\ 
(1,0,1),\ 
(0,1,1).
$$

\noindent
Thus, even though we started with a group involving lots of seventh roots of unity, the resulting cubic surface and braid orbit only involve the integers $0$ and $1$.

This braid orbit turns out to have a very nice realisation from the $G_2$ point of view (which is perhaps not surprising given the strong link to the number $7$):

Suppose we take three lines passing through a single point in the Fano plane (see Figure \ref{Fano Plane}). 
For example the three lines:
$$(e_1,e_3,e_7),\qquad (e_2,e_6,e_7)\qquad (e_4,e_5,e_7)$$
through $e_7$.
Then we obtain three imaginary 
octonions of norm $3$:
$$v_1=e_1+e_3+e_7,\qquad v_2=e_2+e_6+e_7,\qquad v_3=e_4+e_5+e_7$$
and, via Proposition \ref{prop CeqO}, we can consider the corresponding elements of the six dimensional semisimple conjugacy class $\cC\subset G_2$, obtained by conjugating by $a(v_i)=(1+v_i)/2$ respectively:
$$g_1 = T_{a(v_1)},\qquad g_2 = T_{a(v_2)},\qquad g_3 = T_{a(v_3)}\quad\in \cC\subset G_2.$$

Note that the fact that we have chosen three lines passing through a single point implies that the three elements $(1+v_i)/2$ are all {\em octavian integers} (i.e. they lie in the same maximal order, the $7$-integers, cf. \cite{conway-smith} \S9.2).
The octavian integers form a copy of the $E_8$ root lattice (rescaled so roots have norm $1$) and our
three elements $(1+v_i)/2$ are amongst the $240$ units and so correspond to $E_8$ roots.
This implies that with respect to a $\IZ$-basis of the octavian integers the elements $g_1,g_2,g_3$ are represented by matrices with $\IZ$ entries. In fact they generate a finite simple group, and braid to give the Klein orbit of size seven:

\begin{thm}
1) The elements $g_1,g_2,g_3\in G_2(\IC)$ obtained from three lines passing through a single point in the Fano plane, generate a finite subgroup of $G_2(\IC)$ isomorphic to the finite simple group 
$G_2(2)'\cong U_3(3)$ of order $6048$,

2) the braid group orbit of the conjugacy class of 
$(g_1,g_2,g_3)$ in $\cC^3/G=\IM$ is of size seven and lives in a $G_2(\IC)$ character variety isomorphic to the Klein cubic surface.
\end{thm}
\pf
For 1), by construction we obtain a subgroup of the full automorphism group ($\cong G_2(2)$) of the ring of octavian integers. Computing the order shows it is the index two subgroup $\cong G_2(2)'$.
(This was actually our starting point, using Griess's tables \cite{griess-G2} to see that $G_2(2)'$ contains elements of the six dimensional class $\cC\subset G_2(\IC)$.)

2) By construction $\langle v_i,v_j\rangle=1$ 
if $i\ne j$ so that $p_1=p_2=p_3=1$ and 
also we compute $p_4=\langle v_1,v_2\cdot v_3\rangle = -2$.
Thus from \eqref{eq: bdefn} and \eqref{eq: p to xyz}
$b=-1$ and $x=y=z=0$ so that  in turn (by \eqref{eq: symfcs}) $c=0$. Since these parameters are off of all the discriminants we get an isomorphism from the corresponding component of the $G_2$ character variety to the Klein cubic surface. By Corollary \ref{cor: breq} the braid group orbits match up. \epf

Thus the Klein cubic surface \eqref{eq: klein cubic} is related to both the simple group of order $168=2^3\cdot 3\cdot 7$ and 
the simple group of order $6048=2^5\cdot 3^3\cdot 7$. 

\begin{rmk}
If we consider, from the $G_2$ point of view, the symmetric Fricke cubic containing the size $18$ braid orbit of \cite{octa} p.104 (also related to the $237$ triangle group), then we find
$p_1=p_2=-1,p_3=1-4\cos(\pi/7),p_4=1+p_3$. 
Any corresponding triple of elements 
$g_1,g_2,g_3\in \cC\subset G_2$ generate an infinite group:
this value of $p_3$ implies $g_1^2g_2$ has an eigenvalue $x$ with minimal polynomial ${x}^{6}-2\,{x}^{5}+2\,{x}^{4}-3\,{x}^{3}+2\,{x}^{2}-2\,x+1$, and this polynomial has a real root $>1$, so $x$ is not a root of unity.
Thus the speculation/conjecture (of \cite{octa} p.104) that there is a realisation of this Fricke surface relating this braid orbit to a finite group, remains open. 
\end{rmk}

{\small
\noindent{\bf Acknowledgments.}
The first named author is grateful to Nigel Hitchin and Zhiwei Yun for helpful remarks, to the CRM Barcelona for hospitality in May 2012
and was supported by ANR grants 
08-BLAN-0317-01/02, 09-JCJC-0102-01.
}

\renewcommand{\baselinestretch}{1}              %
\normalsize
\bibliographystyle{amsplain}    \label{biby}
\bibliography{../thesis/syr} 

\def\cprime{$'$} \def\cprime{$'$} \def\cprime{$'$} \def\cprime{$'$}
\providecommand{\bysame}{\leavevmode\hbox to3em{\hrulefill}\thinspace}
\providecommand{\MR}{\relax\ifhmode\unskip\space\fi MR }
\providecommand{\MRhref}[2]{%
  \href{http://www.ams.org/mathscinet-getitem?mr=#1}{#2}
}
\providecommand{\href}[2]{#2}
\begin{thebibliography}{10}

\bibitem{AL-p6ims}
D.~Arinkin and S.~Lysenko, \emph{Isomorphisms between moduli spaces of {${\rm
  SL}(2)$}-bundles with connections on {${\bf P}^1\setminus \{x\sb 1,\cdots,
  x\sb 4\}$}}, Math. Res. Lett. \textbf{4} (1997), no.~2-3, 181--190.

\bibitem{AB83}
M.F. Atiyah and R.~Bott, \emph{The {Y}ang-{M}ills equations over {R}iemann
  surfaces}, Phil. Trans. R. Soc. London \textbf{308} (1983), 523--615.

\bibitem{Aud95long}
M.~Audin, \emph{Lectures on gauge theory and integrable systems}, Gauge Theory
  and Symplectic Geometry (J.~Hurtubise and F.~Lalonde, eds.), NATO ASI Series
  C: Maths \& Phys., vol. 488, Kluwer, 1995.

\bibitem{baez-octonions}
J.~C. Baez, \emph{The octonions}, Bull. Amer. Math. Soc. (N.S.) \textbf{39}
  (2002), no.~2, 145--205.

\bibitem{BD96}
A.~Beilinson and V.~Drinfeld, \emph{Quantisation of {H}itchin's fibration and
  {L}anglands program}, Algebraic and Geometric methods in Mathematical Physics
  (A.~Boutet~de Monvel and V.A. Marchenko, eds.), Math. Phys. Studies, vol.~19,
  Kluwer Ac. Pub., Netherlands, 1996.

\bibitem{pecr-nopham}
P.~P. Boalch, \emph{Painlev\'e equations and complex reflections}, Ann. Inst.
  Fourier \textbf{53} (2003), no.~4, 1009--1022, Proceedings of conference in
  honour of Fr\'ed\'eric Pham, Nice, July 2002.

\bibitem{k2p-short}
\bysame, \emph{From {K}lein to {P}ainlev\'e via {F}ourier, {L}aplace and
  {J}imbo}, Proc. LMS \textbf{90} (2005), no.~3, 167--208.

\bibitem{icosa-short}
\bysame, \emph{The fifty-two icosahedral solutions to {P}ainlev\'e {VI}},
  Crelle \textbf{596} (2006), 183--214.

\bibitem{srops-short}
\bysame, \emph{Six results on {P}ainlev\'e {VI}}, S\'emin. Congr., vol.~14,
  SMF, Paris, 2006, pp.~1--20.

\bibitem{octa}
\bysame, \emph{Some explicit solutions to the {R}iemann--{H}ilbert problem},
  IRMA Lectures in Mathematics and Theoretical Physics \textbf{9} (2006),
  85--112, math.DG/0501464.

\bibitem{rsode}
\bysame, \emph{Irregular connections and {K}ac--{M}oody root systems}, 2008,
  arXiv:0806.1050.

\bibitem{logahoric}
\bysame, \emph{Riemann--{H}ilbert for tame complex parahoric connections},
  Transform. Groups \textbf{16} (2011), no.~1, 27--50, arXiv:1003.3177.

\bibitem{ihptalk}
\bysame, \emph{Hyperk\"ahler manifolds and nonabelian {H}odge theory of
  (irregular) curves}, 2012, text of talk at Institut Henri Poincar\'e,
  arXiv:1203.6607.

\bibitem{slims-short}
\bysame, \emph{Simply-laced isomonodromy systems}, Publ. Math. I.H.E.S.
  \textbf{116} (2012), no.~1, 1--68.

\bibitem{conway-smith}
J.~H. Conway and D.~A. Smith, \emph{On quaternions and octonions: their
  geometry, arithmetic, and symmetry}, A K Peters Ltd., Natick, MA, 2003.

\bibitem{CB-Shaw}
W.~Crawley-Boevey and P.~Shaw, \emph{Multiplicative preprojective algebras,
  middle convolution and the {D}eligne-{S}impson problem}, Adv. Math.
  \textbf{201} (2006), no.~1, 180--208.

\bibitem{FrickeKleinI}
R.~Fricke and F.~Klein, \emph{Vorlesungen \"uber die {T}heorie der automorphen
  {F}unktionen. {I}}, Druck und Verlag von B. G. Teubner, Leipzig, 1897, p.
  366.

\bibitem{FulHar}
W.~Fulton and J.~Harris, \emph{Representation theory}, GTM, vol. 129, Springer,
  1991.

\bibitem{goldman-top}
W.~M. Goldman, \emph{Topological components of spaces of representations},
  Invent. Math. \textbf{93} (1988), no.~3, 557--607.

\bibitem{griess-G2}
R.~L. Griess, Jr., \emph{Basic conjugacy theorems for {$G_2$}}, Invent. Math.
  \textbf{121} (1995), no.~2, 257--277.

\bibitem{Hit-tei}
N.~J. Hitchin, \emph{Twistor spaces, {E}instein metrics and isomonodromic
  deformations}, J. Differential Geom. \textbf{42} (1995), no.~1, 30--112.

\bibitem{Hit-G2}
\bysame, \emph{Langlands duality and {$G_2$} spectral curves}, Q. J. Math.
  \textbf{58} (2007), no.~3, 319--344.

\bibitem{Iwas-modular}
K.~Iwasaki, \emph{A modular group action on cubic surfaces and the monodromy of
  the {P}ainlev\'e {VI} equation}, Proc. Japan Acad., Ser. A \textbf{78}
  (2002), 131--135.

\bibitem{kac-1969}
V.~G. Kac, \emph{Automorphisms of finite order of semisimple {L}ie algebras},
  Funkcional. Anal. i Prilo\v zen. \textbf{3} (1969), no.~3, 94--96.

\bibitem{Kron.ale}
P.~B. Kronheimer, \emph{The construction of {ALE} spaces as hyper-{K}\"ahler
  quotients}, J. Differential Geom. \textbf{29} (1989), no.~3, 665--683.

\bibitem{lorenz-mit}
M.~Lorenz, \emph{Multiplicative invariant theory}, Encyclopaedia of
  Mathematical Sciences, vol. 135, Springer-Verlag, Berlin, 2005, Invariant
  Theory and Algebraic Transformation Groups, VI.

\bibitem{Magnus}
W.~Magnus, \emph{Rings of {F}ricke characters and automorphism groups of free
  groups}, Math. Z. (1980), 91--103.

\bibitem{Manin-P6}
Yu.~I. Manin, \emph{Sixth {P}ainlev\'e equation, universal elliptic curve, and
  mirror of {$\bold P^2$}}, Geometry of differential equations, AMS Transl.
  Ser. 2, vol. 186, AMS, Providence, RI, 1998, pp.~131--151.

\bibitem{naruki-cr-looij}
I.~Naruki, \emph{Cross ratio variety as a moduli space of cubic surfaces},
  Proc. London Math. Soc. (3) \textbf{45} (1982), no.~1, 1--30, With an
  appendix by E. Looijenga.

\bibitem{oblomkov-cubics}
A.~Oblomkov, \emph{Double affine {H}ecke algebras of rank 1 and affine cubic
  surfaces}, Int. Math. Res. Not. (2004), no.~18, 877--912.

\bibitem{OkaPVI}
K.~Okamoto, \emph{Studies on the {P}ainlev\'e equations. {I}. {S}ixth
  {P}ainlev\'e equation {$P\sb {{\rm VI}}$}}, Ann. Mat. Pura Appl. (4)
  \textbf{146} (1987), 337--381.

\bibitem{schwarz-regular}
G.~W. Schwarz, \emph{Representations of simple {L}ie groups with regular rings
  of invariants}, Invent. Math. \textbf{49} (1978), no.~2, 167--191.

\bibitem{schwarz-g2bams}
\bysame, \emph{Invariant theory of {$G_{2}$}}, Bull. Amer. Math. Soc. (N.S.)
  \textbf{9} (1983), no.~3, 335--338.

\bibitem{serre-kac-coords}
J.-P. Serre, \emph{Coordin\'ees de {K}ac},
  http://www.college-de-france.fr/site/historique/essai.htm, 2006.

\bibitem{vogt}
H.~Vogt, \emph{Sur les les invariants fondamentaux des \'equations
  diff\'erentielles lin\'eaires du seconde ordre}, Ann Sci. Ecole Norm. Sup.
  \textbf{3} (1889), no.~6, 3--72.

\bibitem{yamakawa-mpa}
D.~Yamakawa, \emph{Geometry of multiplicative preprojective algebra}, Int.
  Math. Res. Pap. IMRP (2008), 77pp.

\end{thebibliography}

\ 

\noindent
D\'epartement de Math\'ematiques, \\
B\^{a}timent 425, \\
Universit\'e Paris-Sud \\
91405 Orsay, France\\

\noindent
philip.boalch@math.u-psud.fr 

\noindent
robert.paluba@math.u-psud.fr

\end{document}